\documentclass[11pt]{amsart}
\usepackage{amsmath,amsfonts,amssymb,amsthm,epsfig,color}
\newcommand{\bb}{\mathbb}

\newcommand{\h}{\bb H}
\newcommand{\Z}{\bb Z}
\newcommand{\R}{\bb R}
\newcommand{\N}{\bb N}

\newcommand{\Q}{\bb Q}
\newcommand{\ff}{\bb F}
\newcommand{\A}{\bb A}

\newcommand{\hh}{\mathcal H}

\newtheorem{Theorem}{Theorem}
\numberwithin{Theorem}{section}
\newtheorem{Cor}[Theorem]{Corollary}
\newtheorem{Prop}[Theorem]{Proposition}
\newtheorem{lemma}[Theorem]{Lemma}

\newtheorem*{lemma*}{Lemma}

\newtheorem*{theorem*}{Theorem}
\newtheorem{Def}{Definition}
\numberwithin{equation}{section}
\begin{document}
\title{Logarithm laws for unipotent flows, I}
\author{J.~S.~Athreya and G.~A.~Margulis}
\subjclass[2000]{primary: 327A17, secondary: 11H16}
\email{jathreya@gmail.com}
\email{margulis@math.yale.edu}
\address{Dept. of Mathematics, Yale University, New Haven, CT}
\thanks{J.S.A. supported by NSF grant DMS
    0603636.}
  \thanks{G.M. supported by NSF grant DMS 0801195} 
\begin{abstract}
We prove analogues of the logarithm laws of Sullivan and Kleinbock-Margulis in the context of unipotent flows. In particular, we obtain results for one-parameter actions on the space of lattices $SL(n, \R)/SL(n, \Z)$. The key lemma for our results says the measure of the set of unimodular lattices in $\R^n$ that does not intersect a  `large' volume subset of $\R^n$  is `small'. This can be considered as a `random' analogue of the classical Minkowski theorem in the geometry of numbers. 
\end{abstract}

\maketitle

\section{Introduction}

An important source of examples for homogeneous dynamics on non-compact spaces is given by the space of unimodular lattices in Euclidean spaces. Let $n \geq 2$, and $X_n = SL(n, \R)/SL(n, \Z)$ denote the space of unimodular lattices in $\R^n$. Let $\mu = \mu_n$ be Haar measure on $X_n$. Define $\alpha_1 : X_n \rightarrow \R^+$ by $$\alpha_1 (\Lambda) := \sup_{ 0 \neq v \in \Lambda} \frac{1}{||v||},$$\noindent where $|| . ||$ denotes Euclidean norm on $\R^n$. By Mahler's compactness criterion, $\alpha_1$ is a \emph{proper} (unbounded off compact sets) function on $X_n$.

Given a (non-compact) one-parameter subgroup $\{g_t\} \in SL(n, \R)$, its acts on $X_n$ by left multiplication. By Moore ergodicity~\cite{Moore}, the action is ergodic with respect to $\mu$, and so almost every trajectory is dense (and in fact uniformly distributed). In particular, we have, for almost every $\Lambda \in X_n$, \begin{equation}\label{escape} \limsup_{t \rightarrow \infty}  \alpha_1(g_t \Lambda) = \infty\end{equation}

A natural question is at what rate the `escape' in equation (\ref{escape}) occurs. For diagonalizable flows, fine properties of cusp excursions were first studied in~\cite{Sullivan} by Sullivan (in the context of finite volume hyperbolic manifolds) and later, in the more general context of the actions of one-parameter diagonalizable subgroups on non-compact finite-volume homogeneous spaces, by Kleinbock-Margulis~\cite{KM}.

In the particular context of the space of lattices, the following theorem was proved in~\cite{KM}:

\begin{Theorem}\label{loglawzero}(\cite{KM}, Theorem 1.7 and Prop 7.1) Fix notation as above. Let $\{a_t\}_{t \in \R}$ denote a diagonalizable (over $\R$) 1-parameter subgroup of $SL(n, \R)$. Then for $\mu$-a.e. $\Lambda \in X_n$, $$\limsup_{t \rightarrow \infty} \frac{\log \alpha_1(a_t \Lambda)}{\log t} = \frac{1}{n}.$$
\end{Theorem}

\vspace{.1in}

Our main result (Theorem~\ref{lattices}) extends this to the context of \emph{unipotent} flows. Our main technical tool (Theorem~\ref{measure} and Corollary~\ref{measurecor}) is a `random' analogue of Minkowski's theorem on lattice points in convex bodies, and is a result of independent interest in the geometry of numbers. We also obtain results on specific sets of orbits in the case $n=2$.

%


%


This paper is the first of a series of two. In the sequel~\cite{AM}, we use more dynamical arguments to obtain results for horospherical actions and one-parameter flows on more general homogeneous spaces $G/\Gamma$. For an example of our results, see subsection~\ref{horoactions}.

\subsection{Organization} This paper is organized as follows: In section~\ref{results}, we state our main results for the space of lattices and geometry of numbers. In section~\ref{upper}, we explore some connections to diophantine approximation (Propositions~\ref{specialdioph} and ~\ref{upperlowerhorodioph}) in the special case $n=2$. In section~\ref{geomnum}, we prove our geometry of numbers result Theorem~\ref{measure} using a formula of Rogers~\cite{Rogers}. In section~\ref{latticeproofs}, we prove our main dynamical results (Theorem~\ref{lattices} and Propostion~\ref{sharp}). Finally, in section~\ref{questions} we describe some interesting questions raised by our work.

\subsection{Acknowledgements:} We would like to thank  Dmitry Kleinbock and Yair Minsky for valuable discussions. Akshay Venkatesh pointed out to us how to extend Theorem~\ref{measure} to the case $n=2$. We would also like to thank Anish Ghosh for pointing us in the direction of the crucial reference~\cite{Rogers}. 

Originally this and ~\cite{AM} were one paper. During the revisions of that earlier manuscript, it was decided to split the paper. The referee's reports on the undivided paper were invaluable in improving the presentation, correcting mistakes, and suggesting a strengthening (Proposition~\ref{upperlowerhorodioph}) of Proposition~\ref{specialdioph}.

\section{Statement of results}\label{results} 


Let $X_n$, $\mu$, and $\alpha_1$ be as above. Let $\{u_t\}_{t \in \R}$ be a unipotent 1-parameter subgroup of $SL(n, \R)$.

\begin{Theorem}\label{lattices} For $\mu$- a.e. $\Lambda \in X_n$, $$\limsup_{t \rightarrow \infty} \frac{\log \alpha_1(u_t \Lambda)}{\log t} = \frac{1}{n}.$$ Moreover, for $n=2$ we have $$\limsup_{t
 \rightarrow \infty} \frac{\log \alpha_1(u_t \Lambda)}{\log t} \geq 1/2$$ for \emph{all} $\Lambda$ such that $\{u_t \Lambda\}_{t \in \R}$ is not periodic. 

\end{Theorem}

\medskip

Our main tool in the proof of Theorem~\ref{lattices} is a `random' analogue of Minkowski's convex body theorem. Recall that Minkowski's theorem states that if $A \subset \R^n$ is a convex, centrally symmetric region with $m(A) > 2^n$ ($m$ is Lebesgue measure on $\R^n$), then for all $\Lambda \in X_n$, there is a non-zero vector $v \in \Lambda \cap A$. 

Without the strong assumptions on the geometry of $A$, the result fails. However, one can ask a probabilistic question: given a set $A$ of large measure, what is the probability that a random lattice (chosen according to Haar measure on $X_n$) does not intersect $A$?

\begin{Theorem}\label{measure} Let $n \geq 2$. There is a constant $C_n$ such that  if $A$ is a measurable set in $\R^n$, with $m(A) >0$, $$\mu(\Lambda \in X_n: \Lambda \cap A = \emptyset) \le \frac{C_n}{m(A)}.$$\end{Theorem}

That is, the measure of the set of lattices that `miss' a set $A$ is bounded above by a quantity inversely proportional to the volume of $A$.

\begin{Cor}\label{measurecor} Let $n \geq 2$. Let $\{A_k\}_{k \in \N}$ be a sequence of measurable sets in $\R^n$, such that $m(A_k) \rightarrow \infty$. Then $$\lim_{k \rightarrow \infty} \mu(\Lambda \in X_n: \Lambda \cap A_k = \emptyset) = 0.$$\end{Cor}

Note that this also yields the fact that for a set $A$ of infinite measure, almost every lattice in $X_n$ will intersect $A$.

In section~\ref{geomnum}, we prove Theorem~\ref{measure} using a key lemma of Rogers~\cite{Rogers} for $n \geq 3$. In subsection~\ref{Eisenstein} we prove Theorem~\ref{measure} in the case $n=2$ following a different approach suggested by A.~Venkatesh. We obtain explicit constants $C_n$, though we do not claim optimality. We then use these results in section~\ref{latticeproofs} to prove Theorem~\ref{lattices}.

\section{Upper bounds and diophantine approximation}\label{upper}

We first explore some connections to diophantine approximation in the special case $\Gamma = SL(2, \Z)$, and the unipotent subgroup \begin{equation}\label{horodef}\left\{h_t =\left(\begin{array}{cc} 1 & t \\ 0 & 1
\end{array}\right): t \in \R\right\}.\end{equation}

\medskip

\noindent\textbf{Remark:} In general, our proofs will be written with the acting unipotent subgroup in Jordan normal form. There is no loss of generality since our results are invariant under conjugation.

Given $s \in \R$, let \begin{equation}\label{slattice}\Lambda_s : =\left(\begin{array}{cc} 1 & 0 \\ s & 1
\end{array}\right)\Z^2\end{equation} \noindent $\Lambda_s$ is a unimodular lattice in $\R^2$, and can be viewed as a point in $X_2 = SL(2, \R)/SL(2, \Z)$. 

Recall the definition of the \emph{approximation exponent} $\mu$ of a real number $s$: \begin{equation}\label{approxexp} \mu = \mu(s) : = \sup \{ \nu \in \R: |s - \frac{p}{q}| < \frac{1}{q^\nu} \mbox{ has infinitely many solutions } \frac{p}{q} \in \Q\}\end{equation} 

\noindent We recall a number of classical facts about $\mu$: 


\begin{enumerate}

\item For almost every $s \in \R$, $\mu(s) = 2$.

\item $\mu(s) = 1$ if and only if $s \in \Q$.

\item For all irrational numbers $s \in \R$, $\mu(s) \geq 2$. 

\item For all algebraic irrational numbers, $\mu(s) =2$ (Roth's theorem).

\item If $\mu(s) = \infty$, we say $\mu$ is a \emph{Liouville} number. This is a set of second category in $\R$.

\end{enumerate}

\begin{Prop}\label{specialdioph} Let $\Lambda_s$ be as above, and $\mu = \mu(s)$. Then for $s \notin \Q$\begin{equation}\label{horodioph}\limsup_{|t| \rightarrow \infty} \frac{\log(\alpha_1(h_t \Lambda_s))}{\log |t|} = 1-\frac{1}{\mu},\end{equation}\noindent and \begin{equation}\label{geoddioph}\limsup_{t \rightarrow \infty} \frac{\log(\alpha_1(g_{-t} \Lambda_s))}{t} = \frac{1}{2}-\frac{1}{\mu}.\end{equation}\end{Prop}

\medskip


\medskip

There is a further refinement of Proposition~\ref{specialdioph}, see Proposition~\ref{upperlowerhorodioph} below.

We also have the following immediate corollary relating to the set of Liouville numbers:

\begin{Cor}\label{liouville} $$\limsup_{|t| \rightarrow \infty} \frac{\log(\alpha_1(h_t \Lambda_s))}{\log |t|} = 1$$ if and only if $s$ is a Liouville number.\end{Cor}

\medskip

\noindent\textbf{Proof of Proposition~\ref{specialdioph}:} 

Let $s \notin \Q$. We first prove the lower bound in (\ref{horodioph}), i.e., that $$\limsup_{|t| \rightarrow \infty}\frac{\log(\alpha_1(h_t \Lambda_s))}{\log |t|} \geq 1-\frac{1}{\mu}.$$ Let $1 < \nu < \mu$. Then there are infinitely many $p, q \in \Z$ with $$|s -\frac{p}{q}| < \frac{1}{|q|^{\nu}}.$$ \noindent Re-writing this equation, we obtain $$|sq - p|^{\frac{1}{1-\nu}} > |q|.$$

\noindent Let $t = t_q = -\frac{q}{sq-p}$. We have $$\alpha_1(h_t \Lambda_s) > \frac{1}{|sq -p|},$$\noindent since $(q, sq-p)^T \in \Lambda_s$ and thus $(0, sq-p)^T \in h_t \Lambda_s$. 

\noindent We also have $$|t| < \frac{|sq-p|^{\frac{1}{1-\nu}}}{|sq-p|} = \frac{1}{|sq-p|^{1 - \frac{1}{1-\nu}}},$$ \noindent since $|sq - p|^{\frac{1}{1-\nu}} > |q|.$ Thus we obtain $$\alpha_1(h_t \Lambda_s) > \frac{1}{|sq -p|} > |t|^{\frac{1}{1 - \frac{1}{1-\nu}}} = |t|^{1- \frac{1}{\nu}},$$ \noindent and thus $$\limsup_{|t| \rightarrow \infty} \frac{ \alpha_1(h_t \Lambda_s)}{|t|^{1-\frac{1}{\nu}}} \geq 1$$ \noindent for all $1 < \nu < \mu$, and thus, we obtain $$\limsup_{|t| \rightarrow \infty} \frac{\log(\alpha_1(h_t \Lambda_s))}{\log |t|} \geq 1-\frac{1}{\mu}.$$

To show the opposite inequality, suppose that $$\limsup_{|t| \rightarrow \infty} \frac{\log(\alpha_1(h_t \Lambda_s))}{\log t} > \beta > 1-\frac{1}{\mu}.$$ \noindent Then, there exists a sequence $\{t_k\}$ with $|t_k | \rightarrow \infty$, and $$\log(\alpha_1(h_{t_k} \Lambda_s) )> \beta \log |t_k|,$$ \noindent Since elements of the lattice $h_t \Lambda_s$ have the form $(m +t(sm+n), sm+n)^T$, this implies that there are $(m_k, n_k) \in \Z^2$ such that $$|m_k + t_k(sm_k +n_k)| < \frac{1}{|t_k|^{\beta}},$$ and \noindent $$|s m_k + n_k| <\frac{1}{|t_k|^{\beta}}.$$

\noindent Now, $$|m_k| < |t_k||sm_k + n_k| + |t_k|^{-\beta},$$ \noindent and $$|t_k||sm_k + n_k| < |t_k|^{1-\beta}.$$

\noindent Let $\nu < 1+ \frac{\beta}{1-\beta}$. We rewrite $|s m_k + n_k| < |t_k|^{-\beta}$ as $$|s m_k + n_k| < |m_k|^{1-\nu} |m_k|^{\nu -1} |t_k|^{-\beta},$$\noindent and thus obtain $$|s m_k + n_k| < |m_k|^{1-\nu} (|t_k|^{1-\beta} + |t_k|^{-\beta})^{\nu-1} |t_k|^{-\beta}.$$

\noindent Since $\nu -1 < \frac{\beta}{1-\beta}$, the product of the second two terms goes to $0$ as $|t_k| \rightarrow \infty$, so for all $\nu < 1+ \frac{\beta}{1-\beta}$, $$|s- p/q| < q^{-\nu}$$ \noindent has infinitely many solutions, i.e. $$1+ \frac{\beta}{1-\beta} \geq \mu(s).$$ \noindent After some simple algebra, this yields $\beta \le 1- \frac{1}{\mu}$, in contradiction to our assumption. 

\medskip\noindent While (\ref{geoddioph}) is well known, we derive a simple proof: let $$\kappa = \limsup_{t \rightarrow \infty} \frac{\log(\alpha_1(g_{-t} \Lambda_s))}{t} .$$ \noindent We first show $\kappa \geq \frac{1}{2} - \frac{1}{\mu}$. Let $2 \le \nu < \mu$ (unless $\mu =2$, then we take $\nu = \mu =2$), so there are infinitely many $p, q \in \Z$ with $$|s -\frac{p}{q}| < \frac{1}{q^{\nu}}.$$

\noindent Since $(q, sq-p)^T \in \Lambda_s$, $(e^{-t/2}q, e^{t/2}(sq-p)) \in g_{-t}\Lambda_s$. To simplify calculations, we consider $\alpha_1$ as defined by the $L^1$-norm on $\R^2$ (this does not make a difference for us since we are interested in $\log \alpha_1$, and any two norms differ by a multiplicative constant), and note that $|sq-p| < q^{1-\nu}$, so that if we minimize \begin{equation}\label{shortgeod}e^{-t/2}q + e^{t/2}q^{1-\nu}\end{equation}\noindent in $t$, we will find a time where $\alpha_1(g_{-t} \Lambda_s)$ is large.

Differentiating (\ref{shortgeod}) and setting the result to $0$, we obtain $t_q= \nu \log q$. So we have $$\alpha_1(g_{-t_q}\Lambda_s) > q^{\frac{\nu}{2}-1} = e^{\frac{t}{\nu} (\frac{\nu}{2} -1)},$$\noindent and thus $$\frac{\log(\alpha_1(g_{-t_q} \Lambda_s))}{t_q} > \frac{\frac{\nu}{2}-1}{\nu} = \frac{1}{2} - \frac{1}{\nu},$$ \noindent so $\kappa \geq \frac{1}{2} - \frac{1}{\nu}$ for all $\nu < \mu$, so $$\kappa \geq \frac{1}{2} - \frac{1}{\mu}.$$ \noindent

To show the opposite inequality, let $\frac{1}{2} \geq \kappa > \lambda$. By assumption, there is a sequence $t_n \rightarrow \infty$ and $(q_n, p_n) \in \Z^2$ such that $(e^{-t_n/2}q_n, e^{t_n/2}(sq_n-p_n))^T$ has length less than $e^{-\lambda t_n}$. In the sequel, we drop the subscript $n$ for brevity. 

Thus, we have $$qe^{-t/2} < e^{-\lambda t}$$ \noindent and $$(qs+p) e^{t/2} < e^{-\lambda t}.$$ \noindent Simplifying, we obtain $q < e^{-(\lambda -\frac{1}{2})t}$ and $(qs+p) < e^{-(\lambda +\frac{1}{2})t}$. Since $\lambda < 1/2$, we have $$q^{\frac{1}{\lambda -\frac{1}{2}}}  > e^{-t} > (qs+p)^{\frac{1}{\lambda + \frac{1}{2}}}.$$\noindent Thus, $$q^{\frac{\lambda+\frac{1}{2}}{\lambda-\frac{1}{2}}} > |qs +p|,$$ \noindent so $$\mu(s) > 1 - \frac{\lambda+\frac{1}{2}}{\lambda-\frac{1}{2}} = \frac{2}{1-2\lambda},$$ which implies $\mu \geq \frac{2}{1-2\kappa}$, and thus, $$\kappa \le \frac{1}{2} - \frac{1}{\mu}$$\noindent as desired.\qed\medskip

It is natural to ask whether one can remove the absolute values from (\ref{horodioph}). Perhaps surprisingly, the answer is no. There is in fact a nice diophantine interpretation of the different limiting behavior of $$\frac{\log \alpha_1(h_t \Lambda_s)}{\log t}$$ as  $t \rightarrow \infty$ and $t \rightarrow -\infty$. 

We define the \emph{upper} and \emph{lower} \emph{approximation exponents} $\mu^{\pm}(s)$ by  

\begin{equation}\label{upperrapproxexp} \mu^{+} = \mu^{+}(s) : = \sup \{ \nu \in \R: 0 < \frac{p}{q} - s < \frac{1}{q^\nu} \mbox{ has infinitely many solutions } \frac{p}{q} \in \Q \}\end{equation} 

\begin{equation}\label{lowerrapproxexp} \mu^{-} = \mu^{-}(s) : = \sup \{ \nu \in \R: 0 < s - \frac{p}{q} < \frac{1}{q^\nu} \mbox{ has infinitely many solutions } \frac{p}{q} \in \Q \}\end{equation} 

We have

\begin{Prop}\label{upperlowerhorodioph} Fix notation as in Proposition~\ref{specialdioph}, and let $\mu^{\pm} = \mu^{\pm}(s)$. Then $$\limsup_{t \rightarrow \infty} \frac{\log \alpha_1(h_t \Lambda_s)}{\log t} = \mu^+,$$\noindent and   $$\limsup_{t \rightarrow \infty} \frac{\log \alpha_1(h_{-t} \Lambda_s)}{\log t} = \mu^-.$$\end{Prop}

\medskip

The proof of this result is a simple modification of the proof of Proposition~\ref{specialdioph}. There, given a `good' approximation $\frac{p}{q}$ to $s$, we produced a time $t_q = - \frac{q}{sq-p}$ for which $\alpha_1(h_{t_q}\Lambda_s)$ is large. Now, it is clear that is $t_q$ is positive or negative according as $s < \frac{p}{q}$ or $s > \frac{p}{q}$. Thus, we have that $t \rightarrow \pm\infty$ yields $\mu^{\pm}$.

To see that $\mu^+$ and $\mu^-$ can be different, we recall their interpretation in terms of the continued fraction expansion of the number $s$. Suppose we have $s = [a_0; a_1, a_2, \ldots, a_n, \ldots]$, with $a_0 = [s] \in \Z$, and $a_i \in \N$. Let $\frac{p_n}{q_n}$ is the $n^{th}$ convergent of $s$, and $r_n : = \frac{\log q_{n+1}}{\log q_n}$. Then using the identity  $$|s - \frac{p_n}{q_n}| < \frac{1}{q_n q_{n+1}} = \frac{1}{q_n^{1+r_n}}$$ (and the fact that if $|s-\frac{p}{q}| < \frac{1}{2q^2}$, $\frac{p}{q} = \frac{p_n}{q_n}$ for some $n$), we can see that $$\mu(s) = 1 + \limsup_{n \rightarrow \infty} r_n.$$

We also have $$\frac{p_{2n}}{q_{2n}} < s < \frac{p_{2n+1}}{q_{2n+1}},$$ and so we can see that $$\mu^+(s) = 1 + \limsup_{n \rightarrow \infty} r_{2n+1},$$ $$\mu^-(s) = 1 + \limsup_{n \rightarrow \infty} r_{2n}.$$ 

Now since, $q_{n+1} = a_n q_n + q_{n-1}$, and we have complete freedom over the $a_n$'s, we see that $r_n$ can have different limiting behavior along its even and odd subsequences, leading to different values for $\mu^+$ and $\mu^-$. 

There is a nice geometric interpretation of Propositions~\ref{specialdioph} and ~\ref{upperlowerhorodioph}. Identifying $SL(2, \R)/SL(2, \Z)$ with the unit tangent bundle to the hyperbolic orbifold $\h^2/SL(2, \Z)$, and consider the standard fundamental domain for $SL(2, \Z)$ on the upper half plane, bounded by the unit circle, and the lines $Re z = \pm 1/2$. The lattice $\Lambda_s$ (restricting to $s \in (-1/2, 1/2)$) can be identified with the point $i + s$ and the upward pointing tangent vector. The geodesic trajectory $\{g_{-t} \Lambda_s\}_{t \geq 0}$ then becomes the line segment from $i+s$ down to $s \in \R$, and the horocycle $\{h_t \Lambda_s\}_{t \in \R}$ is the circle through $i+s$, and tangent to the $x$-axis at $s \in \R$. For $t > 0$, we have the half-circle to the right of the geodesic line, and for $t< 0$, the corresponding half-cricle on the left.

The behavior of cusp excursions of these trajectories is given by the \emph{Ford circles} they pass through. The Ford circles are circles tangent to the $x$-axis at rational points $\frac{p}{q}$ of radius $\frac{1}{2q^2}$. Clearly, the circles the geodesic trajectory passes through correspond to the convergents $\frac{p_n}{q_n}$ associated to $s$, and how far $g_{-t} \Lambda_s$ goes into these circles yields $\mu(s)$. The horocycle on the right ($t>0$) can only encounter approximants to the right of $s$,  yielding $\mu^+(s)$; similarly the left horocycle ($t<0$) yields $\mu^-(s)$.




\section{Geometry of numbers}\label{geomnum}

Let $X_n = SL(n, \R)/SL(n, \Z)$ denote the space of covolume-$1$ lattices in $\R^n$. Let $\mu = \mu_n$ denote Haar measure on $X_n$. Let $B_c(\R^n)$ be the space of measurable, real-valued, bounded, compactly-supported functions on $\R^n$. As is usually done in the geometry of numbers, for $f \in B_c(\R^n)$, define $\widehat{f}: X_n \rightarrow \R$ by $$\widehat{f}(\Lambda) = \sum_{0 \neq v \in \Lambda} f(v).$$

\noindent We will investigate how the $L^2$ norm of $\widehat{f}$ relates to the $L^1$ and $L^2$ norms of $f$. When we consider the $L^1$ norm of $\widehat{f}$, we have the \emph{Siegel Integral Formula:}

\begin{theorem*}(Siegel)$$\int_{X_n} \widehat{f} d\mu = \int_{\R^n} f(v)
dm(v),$$ where $m$ is Lebesgue measure on $\R^n$.\end{theorem*}

The following estimate is an easy corollary of the above result:

\begin{lemma}\label{dl} There is a constant $D_n$ such that for any $r>0$, $$\mu(\Lambda \in X_n: \log\alpha_1(\Lambda) > r) \le D_n e^{-nr}.$$\end{lemma}

\noindent\textbf{Proof:} Apply the Siegel formula to the indicator function of the ball of radius $e^{-r}$\qed\medskip

For $L^2$-norms we have two main results, special cases of results originally due to Rogers~\cite{Rogers}. First, we specialize to the case where $f = I_A$, the indicator function of a measurable set $A$. Define $a := || I_A||_1 = m(A)$. Let $B_a$ be the Euclidean ball around the origin $0$ with $m(B_a) = a$.

\medskip

\noindent\textbf{Remark:} In the rest of this section, we will assume $A$ is bounded. We indicate how to extend to the general case here. Let $A_r : = A \cap B(0, r)$ where $B(0, r)$ is the ball of radius $r$ around the origin. Clearly we have $m(A_r) \rightarrow m(A)$ as $r \rightarrow \infty$ ($\mu(A)$ is possibly infinity). Also, we have $$\mu(\Lambda \in X_n: \Lambda \cap A = \emptyset) \le \mu(\Lambda \in X_n: \Lambda \cap A_r = \emptyset).$$\noindent We have $$\mu(\Lambda \in X_n: \Lambda \cap A_r = \emptyset) \le \frac{C_n}{m(A_r)}$$\noindent by Theorem~\ref{measure}. Thus, for any $r >0$, we have $$\mu(\Lambda \in X_n: \Lambda \cap A = \emptyset) \le \frac{C_n}{m(A_r)},$$\noindent and passing to the limit, we obtain Theorem~\ref{measure} for arbitrary $A$.

\medskip

\begin{lemma}\label{rogers}(\cite{Rogers},Theorem 1) For $n \geq 3$, $$\int_{X_n} \widehat{I_A}^2 d\mu \le \int_{X_n} \widehat{I_{B_a}}^2 d\mu.$$\end{lemma}

\noindent\textbf{Remark}: Rogers works with general measurable functions and their \emph{spherical symmetrizations}. The above result is Theorem 1 of ~\cite{Rogers} specified to indicator functions.

\medskip

This result is a consequence of the following integral formula, also found in~\cite{Rogers}:

\begin{lemma}\label{rogersformula}(\cite{Rogers}, Lemma 1) Let $n \geq 3$. Let $f$ be a non-negative measurable function on $\R^n$. Then $$\int_{X_n} \widehat{f}^2(\Lambda) d\mu = \left(\int_{\R^n} fdm\right)^2 + \sum_{k, q \in \Z: (k, q) =1} \int_{\R^n} f(kx)f(qx)dm(x),$$ where $(k,q)$ is the greatest common divisor of $k$ and $q$. \end{lemma}

In fact, the proof of this result (and its analogs for $p$-norms, $p >2$, also given by Rogers) and the proof of the Siegel integral formula (which can be thought of as such a formula for $p=1$) follow on much the same lines. Define a functional $T_p$ on the space $B_c((\R^n \backslash \{0\})^p)$ of bounded, compactly-supported functions on $(\R^n \backslash \{0\})^p$ by $$T_p(h) = \int_{X_n} \sum_{v_1, \ldots, v_p \in \Lambda \backslash {0}} h(v_1, \ldots, v_p)d\mu.$$\noindent This is  $SL(n, \R)$ invariant, and thus the measure that defines it must be the combination of $SL(n, \R)$-invariant measures on $SL(n, \R)$-orbits in $(\R^n \backslash \{0\})^p$. Performing this decomposition and applying it to functions $h \in  B_c(\R^n \backslash \{0\})^p$ defined by $h(v_1, \ldots, v_p) = f(v_1)f(v_2)\ldots f(v_p)$, where $f \in B_c(\R^n\backslash{0})$. 

For $p=1$ there is only one orbit, all of $\R^n\backslash\{0\}$. For $p=2$, there are two types of orbits, one consisting of linearly independent pairs of vectors, which yields the first term ($||f||_1^2 = ||\widehat{f}||_1^2$), and the summation captures the contribution of pairs of linearly dependent vectors, with each term associated to a fixed ratio. For $p>2$, the formula becomes much more complicated. It should be noted that the $p$-norm formula only works for dimensions $n>p$. 

\subsection{Proof of Theorem~\ref{measure}} Let $n \geq 3$. Given a measurable set $A \subset \R^n$, let $\sigma_A = \mu(\Sigma_A)$, where $\Sigma_A = \{\Lambda \in X_n: \Lambda \cap A = \emptyset\}$. We want to calculate an upper bound for $\sigma_A$ in terms of $a:= m(A)$.

We define $g_A: X_n \rightarrow \R$ by $g_A = I_{\Sigma_A^c}$. Let $f_A = \widehat{I_A}$. Note that $f_A = g_A f_A$, because $f_A = 0$, when $g_A =0$. Thus, by Cauchy-Schwarz:

$$\left(\int_{X_n} f_A d\mu\right)^2 \le \left(\int_{X_n} f_A^2 d\mu\right)\left(\int_{X_n} g_A^2 d\mu\right).$$

\noindent The Siegel formula yields $||f_A||_1 = m(A)$, and by Lemma~\ref{rogers},$$\int_{X_n} f_A^2 d\mu \le \int_{X_n} \widehat{I_{B_a}}^2 d\mu.$$\noindent Also, $\sigma_A = 1 - ||g_A||_1 = 1- ||g_A||_2^2$ (since $g_A$ is an indicator function, $g_A^2 = g_A$).  So we have:

$$a^2 \le || \widehat{I_{B_a}}||_2^2 \left(1 - \sigma_A\right),$$

\noindent and thus, $$\sigma_A \le 1 - \frac{a^2}{|| \widehat{I_{B_a}}||_2^2}.$$\noindent

Note that $m(A) = ||\widehat{I_{B_a}}||_1$, so to complete the proof, it suffices to calculate $||\widehat{I_{B_a}}||_2^2$

$$|| \widehat{I_{B_a}}||_2^2 = a^2 + \sum_{k,q \in \Z: (k,q) =1}\int_{\R^n} I_{B_a}(kv) I_{B_a}(kq) dm(v),$$\noindent A simple calculation shows that the sum can be re-written as $$a \sum_{k,q \in \Z: (k,q) =1} \frac{1}{\max(|k|,|q|)^n}.$$\noindent Now, $$\sum_{k,q \in \Z: (k,q) =1} \frac{1}{\max(|k|,|q|)^n} = 8\sum_{k=1}^{\infty} \frac{\varphi(k)}{k^n},$$\noindent where $\varphi$ is the Euler totient function. Since $n \geq 3$, the right-hand sum is convergent. Thus, we obtain that $$|| \widehat{I_{B_a}}||_2^2 = a^2 + C_n a,$$\noindent for $C_n =8\sum_{k=1}^{\infty} \frac{\varphi(k)}{k^n} = 8\frac{\zeta(n-1)}{\zeta(n)}$, and so, $$\sigma_A \le \frac{C_n a}{a^2 + C_n a} \le \frac{C_n}{a}$$ which completes the proof. \qed

\subsection{The $n=2$ case}\label{Eisenstein}

Note that the above proof does not work in the case $n=2$, since Lemma~\ref{rogersformula} does not apply. However, one can still obtain estimates on the $L^2$ norm of an appropriate transform in order to estimate the measure of the set of lattices that do not intersect a fixed set. 

We follow the presentation in Chapter 13 of Lang's book~\cite{Lang} (although, to stay consistent with the notation of our paper, we quotient on the right). Let $G = SL(2, \R)$, $N = \{h_t: t \in \R\}$, $\Gamma = SL(2, \Z)$. Note $G/N = \R^2\backslash \{0\}$.

We define the $\Theta$-transform $Tf$ for any function $\phi \in L^1(\R^2)$ by $$T\phi(\Lambda) = \sum_{v \in \Lambda_{prim}} \phi( v),$$\noindent where $\Lambda \in X_2 = SL(2, \R)/SL(2, \Z)$ is a lattice, and $\Lambda_{prim} \subset \Lambda$ denotes the set of primitive vectors. Let $\mu$ denote the non-normalized Haar measure on $X_2$ (i.e., $\mu (X_2) = \zeta(2)$). Then we have the following modification of the Siegel formula:

\begin{theorem*}(Siegel)$$\int_{X_2} (T\phi) d\mu = \int_{\R^n} \phi$$\end{theorem*}

\vspace{.1in}

\noindent We can also define a map $T^0$ going the other way: given a function $f$ on $G/\Gamma$, we define a function $T_0 f$ on $G/\Gamma$. We view $f$ as a $\Gamma$-periodic function on $G$, and define a function on $G/N$ by constructing an $N$-periodic function on $G$: $$(T^0 f)(g) = \int_{N/\Gamma \cap N} f(gn) dn,$$\noindent where the integral is taken with respect to Haar measure on $N$, normalized to be a probability measure on $N/\Gamma \cap N$.

$T_0$ and $T$ satisfy the following adjointness relation:

\begin{Theorem}\label{adjoint} (\cite{Lang}, Chapter XIII, section 1, Theorem 1) let $\phi \in C^{\infty}_c(G/N)$, $f \in L^2(G/\Gamma)$. Then $$\langle T\phi , f \rangle_{G/\Gamma} = \langle \phi, T^0 f \rangle_{G/N},$$ \noindent where $<,>$ denotes the standard $L^2$-inner product. \end{Theorem}

We also have the following relation:

\begin{Theorem}\label{L2 isometry} (\cite{Lang}, Chapter XIII, section 8, last Corollary) Let $\phi, \psi \in C^{\infty}_c(G/N)$ be even, mean zero functions. Then $$\langle T\phi, T\psi \rangle_{G/\Gamma} =  \langle T^0 T \phi, T^0 T \psi \rangle_{G/N}.$$\end{Theorem}

Combining these results we obtain:

\begin{lemma}\label{L2norm} Let $\phi$ be a mean-zero compactly supported function on $\R^2$. Then 

\begin{equation}\label{L2normeq}||T \phi||_2 \le ||\phi||_2.\end{equation}

\end{lemma}

\noindent\textbf{Proof:} Note that we can drop the requirement of evenness since we can split $\phi = \phi_{even} + \phi_{odd}$, and $T\phi = T\phi_{even}$. 

Now, assuming $\phi \in C^{\infty}_c(G/N)$, equation~\ref{L2normeq} follows from $$||T\phi||^2_2 = \langle T\phi, T\phi \rangle = \langle \phi, T^0T\phi\rangle  \le ||\phi||_2 ||T^0 T\phi||_2 = ||\phi||_2 ||T\phi||_2.$$

To extend this to aribtrary compactly supported $\phi$, let $\{\phi_n\}_{n=0}^{\infty} \subset C_c^{\infty}(G/N)$ be a sequence of even mean-zero functions such that $\phi_n \rightarrow \phi$ in $L^1$ and $L^2$ (this can be done by density of $C_c^{\infty}$ functions). $\{T\phi_n\}_{n=0}^{\infty}$ is a Cauchy sequence, so by passing to a subsequence, we can obtain an almost everywhere convergent subsequence. Now, we have for all $n$, $$||T\phi_n||_2 \le ||\phi_n||_2,$$\noindent and applying the Siegel formula and Fatou's lemma, we obtain $$||T\phi||_2 \le ||\phi||_2.$$\qed\medskip

Now, let $A \subset \R^2$ be measurable. Define $a$ and $B_a$ as above. Let $T_A = T I_{A}$ and $T_B = T I_{B_a}$. $(I_A - I_B)$ is a mean-zero function, so applying Lemma~\ref{L2norm} we have that \begin{equation}\label{AB}||T_A - T_B||_2 \le ||I_A - I_B||_2 \le 2 \sqrt{a}.\end{equation}

Randol (\cite{Randol}, end of section 1) showed that  $$||I_B - a/\zeta(2)||_2^2 = a/\zeta(2)+ O(a/\log a).$$\noindent Thus, for all $a$ sufficiently large, we have that \begin{equation}\label{Bconst}||I_B - a/\zeta(2)||_2 \le 2\sqrt{a}.\end{equation}

Combining equations~\ref{Bconst} and~\ref{AB}, we have \begin{equation}\label{Aconst}||T_A - a/\zeta(2)||_2 \le 4\sqrt{a}.\end{equation}

Letting $p_A = \mu(\Lambda: T_A(\Lambda) = 0)$, we have that $$16a \geq ||T_A - a/\zeta(2)||_2^2 \geq p_A (\frac{a}{\zeta(2)})^2,$$\noindent and thus, $$p_A \le 16\zeta(2)^2/a$$ \noindent for all $a>>0$. This yields Theorem~\ref{measure} for the case $n=2$.\qed\medskip

\section{Lattice results}\label{latticeproofs}

\noindent In this section, we use the results of section~\ref{geomnum} in order to prove Theorem~\ref{lattices}. We recall the statement of Theorem~\ref{lattices}:

\begin{theorem*} For almost every $\Lambda \in X_n$, $$\limsup_{t \rightarrow \infty} \frac{\log \alpha_1(u_t \Lambda)}{\log t} = \frac{1}{n}.$$ Moreover, for $n=2$ and $u_t = h_t$, where $h_t$ is as defined in equation~\ref{horodef}, we have $$\limsup_{t \rightarrow \infty} \frac{\log \alpha_1(u_t \Lambda)}{\log t} \geq 1/2$$ for \emph{all} $\Lambda$ such that $\{u_t \Lambda\}_{t \geq 0}$ is not periodic. 

\end{theorem*}

The proof naturally splits into an upper bound and a lower bound, and the lower bound further splits into the cases $n =2$ and $n>2$. For the rest of this section, we fix notation as in section~\ref{geomnum}.

\subsection{Upper bound:}\label{ub}

\begin{lemma}\label{upperbound} For $\mu$-almost every $\Lambda \in X_n$, $$\limsup_{t \rightarrow \infty} \frac{\log \alpha_1(u_t \Lambda)}{\log t} \le \frac{1}{n}.$$\end{lemma}

\noindent\textbf{Proof:} Let $\epsilon >0$, and $r_k = (\frac{1}{n} + \epsilon)\log k$. We have $$\mu(\Lambda \in X_n: \log \alpha_1(u_k\Lambda) > r_k) \le C_n \frac{1}{k^{1+n\epsilon}},$$\noindent by Lemma~\ref{dl} and the fact that $u_t$ is measure-preserving. But $C_n \frac{1}{k^{1+n\epsilon}}$ is summable in $k$, so by the convergence half of the Borel-Cantelli lemma, for almost every $\Lambda$, $\log \alpha_1(u_k\Lambda) > r_k$ only finitely often. This yields that $$\limsup_{t \rightarrow \infty} \frac{\log \alpha_1(u_t \Lambda)}{\log t} < \frac{1}{n} + \epsilon$$ for almost every $\Lambda$, and since $\epsilon >0$ was arbitary, we obtain our desired upper bound.\qed \medskip

\subsection{Lower bound, $n=2$:}\label{lb2}

\noindent The lower bound can be derived from using the methods described in subsection~\ref{nb3}. However, as a simple example of our methods using the geometry of numbers, we prove the following:

\begin{lemma}\label{lowerbound2} For all $\Lambda$ such that $h_t \Lambda$ is non-periodic  $$\limsup_{|t| \rightarrow \infty} \frac{\alpha_1(h_t \Lambda)}{|t|^{1/2}} \geq 1.$$\end{lemma}

\noindent\textbf{Remark:} It suffices to consider $h_t$ here, since any unipotent one parameter subgroup can be conjugated to $h_t$. 

\medskip

\noindent\textbf{Proof:} Let $\Lambda \in X_2$ be a lattice without a horizontal vector (i.e., it is not $h_t$-periodic). For $k \in \N$, define $$A_k := \{(x,y) \in \R^2: |x| \le \sqrt{k}, |y| \le 1/\sqrt{k}\}.$$\noindent $A_k$ is a convex, centrally symmetric region of area $4$, so by Minkowski's theorem, there is a non-zero point $(x_k, y_k) \in \Lambda$, and moreover, since $\Lambda$ has no horizontal vectors, $y_k \neq 0$ for all $k$. Also, the set of all intersection points $\bigcup_{k \in \N} A_k \cap \Lambda$ is unbounded: otherwise, our lattice $\Lambda$ would have an accumulation point. Thus, by passing to a subsequence if necessary, we can assume that we can pick $\{(x_k, y_k)\}_{k \in \N}$ so that $|x_i| >| x_j|$ and $|y_i| < |y_j|$ for $i >  j$. Let $t_k = -x_k/y_k$. 

$$\alpha_1(h_{t_k} \Lambda) \geq 1/|y_k|,$$

\noindent since $h_{t_k} (x_k, y_k)^T = (0, y_k)^T$. Now $|t_k| = |x_k/y_k| \le \sqrt{k}/|y_k| \le 1/y_k^2$, so $$ \alpha_1(h_{t_k}) \geq 1/|y_k| \geq |t_k|^{1/2},$$ \noindent so we have produced an infinite sequence of times $t_k$ where we achieve our lower bound.\qed

\subsection{Lower bound, $n \geq 3$:}\label{nb3}

\begin{lemma}\label{lowerbound3} For $\mu$-almost every $\Lambda \in X_n$, $$\limsup_{t \rightarrow \infty} \frac{\log \alpha_1(u_t \Lambda)}{\log t} \geq \frac{1}{n}.$$
\end{lemma}

\noindent\textbf{Proof:} As all the essential ideas are contained in the case $n=3$, we describe it first. First, suppose our one-parameter unipotent flow $u_t$ is simply a copy of $h_t$, e.g., $$u_t =\left(\begin{array}{ccc} 1 & 0  &0 \\ 0 & 1 & t\\ 0 & 0 &1
\end{array}\right).$$\noindent Fix $\delta>0$, and define the region $$A_k = \{(x,y, z) \in \R^3: y/z <0,  |x| \le |k|^{-1/3}, |y| \le k^{2/3}, k^{-1/3 - \epsilon} \le |z| \le k^{-1/3+\epsilon}\},$$ where $\epsilon = \epsilon(\delta)$ will be determined later. Now $A_k$ is no longer convex, but a simple volume calculation shows that $m(A_k) = 4 (k^{\epsilon} - k^{-\epsilon})$, and thus we can apply Corollary~\ref{measurecor} to see that $\mu(\Lambda:\Lambda \cap A_k = \emptyset ) \rightarrow 0$. Thus, for almost every $\Lambda$, by passing to a subsequence if needed, we can produce a sequence of distinct non-zero points $\{(x_k, y_k, z_k) \in\Lambda\cap A_k\}$, and with $z_k \neq 0$. Set $t_k = -y_k/z_k$, once again, passing to a subsequence if needed, we can take $t_k \uparrow +\infty$. Now

$$\alpha_1(h_{t_k} \Lambda) \geq \frac{1}{|x_k| + |z_k|} \geq \frac{1}{k^{-1/3}+k^{-1/3+\epsilon}} \geq \frac{1}{2k^{-1/3 +\epsilon}},$$

\noindent  and $$t_k = -y_k/z_k \le k^{2/3}/|z_k|  \le k^{1+\epsilon}.$$

\noindent Thus $$\frac{\log\alpha_1(h_{t_k} \Lambda)}{\log t_k} \geq  \frac{ (1/3-\epsilon) \log k - \log 2}{(1+\epsilon)\log k} = \frac{1}{3} \left(\frac{(1-3\epsilon)}{1+\epsilon} - \frac{\log 2}{(3+3\epsilon)\log k}\right).$$

\noindent Pick $\epsilon = \epsilon(\delta)$, $k_0 = k_0 (\epsilon, \delta)$ so that $$ \frac{1}{3} \left(\frac{(1-3\epsilon)}{1+\epsilon} - \frac{\log 2}{(3+3\epsilon)\log k}\right) > 1/3 - \delta,$$\noindent for all $k > k_0$. 

\noindent We have produced, for each $\delta>0$, a set of full measure $X_{n, \delta} \subset X_n$ such that $$\limsup_{t \rightarrow \infty} \frac{\log \alpha_1(u_t \Lambda)}{\log t} > \frac{1}{3} - \delta.$$\noindent The set $\bigcap_{j \in \N} X_{n, 1/j}$ is a set of full measure on which we have our lower bound.

The other case is when $u_t$ is regular, e.g. 

$$u_t =\left(\begin{array}{ccc} 1 & t  & \frac{t^2}{2}\\ 0 & 1 & t\\ 0 & 0 &1
\end{array}\right) .$$

\noindent Our proof follows on similar lines. Again, fix $\delta>0$, and define the region $$A_k = \{(x,y, z) \in \R^3: y/z < 0, |x- y^2/2z| \le |k|^{-1/3}, |y| \le k^{2/3}, k^{-1/3 - \epsilon} \le |z| \le k^{-1/3+\epsilon}\},$$ where $\epsilon = \epsilon(\delta)$ will be determined later. The reason for this modification is that $$u_t\left(\begin{array}{c}x\\ y \\  z  \end{array}\right) = \left(\begin{array}{c}x + ty +\frac{t^2}{2}z\\ y+tz \\ z \end{array}\right),$$\noindent and if we set $t = -y/z$, we obtain $$\left(\begin{array}{c}x - y^2/2z\\ 0 \\ z \end{array}\right).$$

\noindent The remainder of the calculation proceeds on exactly the same lines as the case $$u_t =\left(\begin{array}{ccc} 1 & 0  &0 \\ 0 & 1 & t\\ 0 & 0 &1
\end{array}\right).$$\noindent The volume calculation is identical, and once we note that $$\alpha_1(h_{t_k} \Lambda) \geq \frac{1}{|x_k - \frac{y_k^2}{2z_k}| + |z_k|} \geq \frac{1}{k^{-1/3}+k^{-1/3+\epsilon}} \geq \frac{1}{2k^{-1/3 +\epsilon}},$$\noindent we can proceed verbatim as above.

For general $n$, the proof follows on much the same lines. Given $u_t$, we decompose it into Jordan blocks, each of the form $$\left(\begin{array}{cccccc}1 & t & t^2/2 & t^3/6 & \ldots & t^k/k!\\ 0 & 1 & t & t^2/2 & \ldots & t^{(k-1)}/(k-1)! \\ 0 & \ldots & . & \ldots & \ldots &\ldots \\ \ldots & 0 & \ldots & . & \ldots & \ldots \\ \ldots & \ldots & 0 & \ldots &1 & t \\ \ldots & \ldots & \ldots & 0& \ldots & 1 \end{array}\right).$$

\noindent Without loss of generality, we assume that the bottom right hand corner of $u_t$ looks like $h_t$. That is, if we apply $u_t$ to a vector $x = (x^{(1)}, x^{(2)} \ldots, x^{(n)})^T$, and ignore the first $(n-2)$ coordinates, we obtain $$u_t x = \left(\begin{array}{c} . \\ . \\ . \\ . \\ x^{(n-1)} + tx^{(n)} \\ x^{(n)}\end{array}\right).$$ \noindent We are not specifying what happens to the other coordinates, but they will be polynomials in $t$ with the $x_i$'s as coefficients, and the constant term in the $i^{th}$ coordinate will be $x_i$. If we set $t_0 = t_0(x) = -x^{(n-1)}/x^{(n)}$, then the $i^{th}$ coordinate can be expressed as the difference of $x^{(i)}$ and a rational function in the other $(n-1)$ variables, call this function $f_i(x)$ (only $x^{(n)}$ will appear in the denominator). That is, we have $$u_{t_0}x = \left(\begin{array}{c} x^{(1)} - f_1(x) \\ x^{(2)} - f_2(x)\\ \ldots \\ x^{(i)} - f_i(x)\\  \ldots \\ 0 \\ x^{(n)} \end{array}\right).$$

\noindent Again, fix $\delta>0$, and let $A_k$ be the set of $x \in \R^n$ so that: 

\begin{enumerate}
\item $x^{(n-1)}/x^{(n)} < 0$
\item $|x^{(i)} - f_i(x)| \le k^{-1/n} \mbox{ for } i < (n-1)$
\item $ |x^{(n-1)}| \le k^{n-1/n}$
\item $|x^{(n)}| \in [k^{-1/n - \epsilon}, k^{-1/n + \epsilon}].$

\end{enumerate}

\noindent where $\epsilon = \epsilon(\delta)>0$ is to be determined later. $m(A_k) = 2^{n-1} (k^{\epsilon} - k^{-\epsilon})$, so $m(A_k) \rightarrow \infty$. By Corollary~\ref{measurecor}, for almost every lattice $\Lambda$, we can produce a sequence of non-zero points $y_k \in (A_k \cap \Lambda)$. Let $t_k = t_0 (y_k)$. Passing to a subsequence if necessary, we can assume $t_k$ is well-defined and positive. The rest of the proof follows along exactly the same lines as above, with $n$ in place of $3$, and $\log(n-1)$ in place of $\log 2$. \qed \medskip

\noindent\textbf{Proof of Theorem~\ref{lattices}:} Combine Lemmas~\ref{upperbound}, ~\ref{lowerbound2}, and ~\ref{lowerbound3}.\qed

\subsection{A sharpening:}\label{sharpen}

\noindent We can refine the above proof on lower bounds to yield the following sharper result:

\begin{Prop}\label{sharp} Fix $n \geq 2$. Let $r: \R \rightarrow \R^+$ be such that there are constants $c, \delta>0$ so that 
\begin{enumerate} 
\item $l(t) := \frac{r(t)}{t^{1/n}}$ is non-decreasing.
\item $\lim_{t \rightarrow \infty} \frac{r(t)}{t^{\frac{1+\delta}{n}}} = 0.$
\item $\liminf_{t \rightarrow \infty} \frac{r(t)t^{2\delta/n}}{r(t^{1+2\delta})} >c.$
\end{enumerate}
Then for almost all $\Lambda \in X_n$, there is a sequence $\{t_k\}_{k \in \N}$ with $t_k \rightarrow + \infty$ such that $$\alpha_1(u_{t_k} \Lambda) \geq \frac{c}{n-1} r(t_k)t_k^{-2\delta/n}.$$
\end{Prop}

\noindent\textbf{Remark:} If $r(t) = t^{1/n}$, we have $c=1$ and $\delta$ arbitrary, and our result essentially reduces to Lemma~\ref{lowerbound3}.

\medskip

\noindent\textbf{Proof:} Again, assume that the bottom right hand corner of $u_t$ looks like $h_t$, and thus there are rational functions $f_i: \R \rightarrow\R$ such that if we have $x \in \R^n$ with $x^{(n)} \neq 0$, and we set  $t_0 = t_0(x) = -x^{(n-1)}/x^{(n)}$, we have $$u_{t_0}x = \left(\begin{array}{c} x^{(1)} - f_1(x) \\ x^{(2)} - f_2(x)\\ \ldots \\ x^{(i)} - f_i(x)\\  \ldots \\ 0 \\ x^{(n)} \end{array}\right).$$

\noindent Let $A_k$ be the set of $x \in \R^n$ so that: 

\begin{enumerate}

\item $x^{(n-1)}/x^{(n)} < 0$
\item $|x^{(i)} - f_i(x)| \le k^{-1/n}l^{-1}(k) = r^{-1}(k) \mbox{ for } i < (n-1)$ 
\item $ k^{\frac{n-1}{n}}l^{-1}(k) \le |x^{(n-1)}| \le k^{\frac{n-1}{n}+ \delta}l^{-1}(k)$
\item $k^{-\frac{1}{n}- \delta}l^{-1}(k) \le |x^{(n)}| \le r^{-1}(k).$

\end{enumerate}

An easy calculation shows that $m(A_k) = 2^nl^{-n}(k)(k^{\delta} + k^{-\delta} -2)$. Now, by condition 2 on $r(t)$, we have that $$\lim_{t \rightarrow \infty} t^{\delta}l^{-n}(t) = \infty,$$\noindent so $m(A_k) \rightarrow \infty$. Applying Corollary~\ref{measurecor}, we obtain that for almost every $\Lambda \in X_n$ there is a sequence $\{x_k\}$ such that $x_k \in A_k \cap \Lambda$. Passing to a subsequence if required, we can assume $t_k = t_0 (x_k) >0$. Note that $$\frac{k^{\frac{n-1}{n}}l(k)}{k^{-1/n} l(k)} \le t_k \le \frac{k^{\frac{n-1}{n}+\delta}l(k)}{k^{-\frac{1}{n} +\delta} l(k)}$$\noindent and thus $$k \le t_k \le k^{1+2\delta}.$$

\noindent Now, $$\alpha(h_{t_k} \Lambda) \geq \frac{1}{|x^{(n)}_k| + (n-2)r^{-1}(k)} \geq \frac{1}{(n-1)} r(k).$$ \noindent Rewrite $$r(k) = r(t_k) t_k^{-2\delta/n} \frac{r(k)t_k^{2\delta/n}}{r(t_k)}.$$\noindent By condition 3 on $r$, for all $k$ sufficiently large $$ \frac{r(k)t_k^{2\delta/n}}{r(t_k)} > \frac{r(k) k^{2\delta/n}}{r(k^{1+2\delta})} > c.$$\qed \medskip

\section{Further questions}\label{questions}

To our knowledge, this paper contains the first results on statistical behavior of the excursions of unipotent flows on non-compact homogeneous spaces. We hope that it will inspire further results in this direction. Some interesting classes of questions include:

\subsection{Horospherical actions}\label{horoactions} In the sequel to this paper~\cite{AM}, we prove results for the excursions of expanding translates of horospherical subgroups for general homogeneous spaces. A nice example of our results is in the context of hyperbolic surfaces:

Let $G = SL(2, \R)$ and let $h_s$ be as in equation~\ref{horodef}. Let $\Gamma \subset SL(2,\R)$ be a non-uniform lattice . Let  $d$ denote distance on the hyperbolic surface $S = \h^2/\Gamma$ ($\h^2$ denotes the upper-half plane with constant curvature $-1$), and $p: M \rightarrow S$ be the natural projection from $M =SL(2, \R) /\Gamma$.

\begin{Theorem}\label{horosl2} Let $H = \{h_t\}_{t \in \R}$.  Fix $y \in S$.Then for all $x \in S$, almost all $\tilde{x} \in p^{-1}(x)$, $$\limsup_{t \rightarrow \infty} \frac{d(p(h_t \tilde{x}), y)}{\log t} = 1.$$ Moreover, for all $\tilde{x} \in M$ such that $H\tilde{x}$ is not closed, $$\limsup_{t \rightarrow \infty} \frac{d(p(h_t \tilde{x}), y)}{\log t} \geq 1.$$ \end{Theorem}

\subsection{Flows on moduli spaces of surfaces and geometry of saddle connections:} Theorem~\ref{AM}, a logarithm law for the horocycle flow on a stratum of abelian differentials, is proved in~\cite{Athreya-Minsky}.

\begin{Theorem}\label{AM} (Athreya-Minsky) Let $\hh$ denote a stratum of the space of abelian differentials on a surface of genus $g >1$. Let $\lambda: \hh \rightarrow \R^+$ be defined by $$\lambda(\omega) = \sup_{v \in V_{sc}(\omega)} \frac{1}{||v||},$$ where $V_{sc}(\omega)$ is the set of (holonomy vectors of) saddle connections on $\omega$. Then for almost every (with respect to the Lebesgue measure) $\omega \in \hh$, $$\limsup_{t \rightarrow \infty} \frac{\log \lambda(h_t \omega)}{\log t} = 1/2,$$ where $h_t$ is as in equation~\ref{horodef}.\end{Theorem}
 
 The main tool is a version of the Minkowski theorem for the set of holonomy vectors of saddle connections. It would be interesting to see if one can get a version of Theorem~\ref{measure} for these sets. 

There is a logarithm law (more analagous to Theorem~\ref{loglawzero}) for the Teichmuller geodesic flow, due to Masur~\cite{Masur}.

\subsection{General shrinking target properties:} A natural complement to the question of excursions into the cusp is the more general question of \emph{shrinking target properties} (or STP's):

Let $(X, \mu)$ be a probability space, and $G$ be a group acting on $X$ by measure-preserving transformations. Let $\mathcal{A}$ be a family of measurable subsets of $X$. Let $\{g_n\}_{n \in \N} \subset G$ be a sequence of group elements.

\begin{Def}\label{STP} We say that $\mathcal{A}$ satisfies the STP for the sequence $\{g_n\}_{n=0}^{\infty}$ if for any sequence $\{A_n\}_{n \in \N} \subset \mathcal{A}$ \begin{equation}
\mu(\{x \in X : g_n x \in A_n \mbox{ infinitely often} \}) =
\left\{ \begin{array}{ll}  1 & \sum_{n=0}^{\infty} \mu(A_n) = \infty \\ 0 & \mbox{otherwise}\end{array}\right.\nonumber\end{equation}\end{Def}

There are many results on STP's for hyperbolic (\cite{CK, Dol,KM, Mau}) and elliptic (\cite{Fayad, Tseng}) systems, but it would be interesting to obtain more results for \emph{parabolic} systems like unipotent flows. 

\subsection{Algebraic groups:} One can consider shrinking target and logarithm law question for $G/\Gamma$ where $G$ is a group over a field $k$. $k$ here could be $\Q_p$, the adeles $\A$, or a field of positive characteristic $\ff_q((t^{-1}))$. These will also naturally relate to the geometry of buildings and trees. Results in this direction for Cartan actions have been obtained in~\cite{AGP, HP}.

\end{document}